\documentclass[12pt]{article}
\usepackage{latexsym,amsfonts,amssymb}
\usepackage[dvips]{color}
\usepackage{colordvi,multicol}

\def \cal{\mathcal}

\newtheorem{thm}{Theorem}[section]
\newtheorem{cor}[thm]{Corollary}

\newtheorem{pro}[thm]{Proposition}

\newtheorem{rem}[thm]{Remark}

\begin{document}
\title{\bf  A dynamic model for the two-parameter\\ Dirichlet process}
\author{}

\maketitle
\date{}

 \centerline{Shui Feng} \centerline{\small Department of Mathematics and Statistics, McMaster
 University,}
\centerline{\small Hamilton, L8S 4K1, Canada} \centerline{\small
E-mail: shuifeng@mcmaster.ca} \vskip 0.7cm \centerline{Wei Sun}
\centerline{\small Department of Mathematics and Statistics,
Concordia University,} \centerline{\small Montreal, H3G 1M8,
Canada} \centerline{\small E-mail: wei.sun@concordia.ca}

\vskip 1cm


\vskip 0.5cm \noindent{\bf Abstract}\quad Let $\alpha=1/2$,
$\theta>-1/2$, and $\nu_0$ be a probability measure on a type
space $S$.  In this paper, we investigate the stochastic dynamic
model for the two-parameter Dirichlet process
$\Pi_{\alpha,\theta,\nu_0}$. If $S=\mathbb{N}$, we show that the
bilinear form
\begin{eqnarray*}
\left\{
\begin{array}{l}
{\cal E}(F,G)=\frac{1}{2}\int_{{\cal P}_1(\mathbb{N})}\langle
\nabla F(\mu),\nabla G(\mu)\rangle_{\mu}
\Pi_{\alpha,\theta,\nu_0}(d\mu),\ \ F,G\in
{\cal F},\\
{\cal F}=\{F(\mu)=f(\mu(1),\dots,\mu(d)):f\in
C^{\infty}(\mathbb{R}^d), d\ge 1\}
\end{array}
\right.
\end{eqnarray*}
is closable on $L^2({\cal
P}_1(\mathbb{N});\Pi_{\alpha,\theta,\nu_0})$ and its closure
$({\cal E}, D({\cal E}))$ is a quasi-regular Dirichlet form. Hence
$({\cal E}, D({\cal E}))$ is associated with a
 diffusion process in ${\cal P}_1(\mathbb{N})$ which is time-reversible with
the stationary distribution $\Pi_{\alpha,\theta,\nu_0}$. If $S$ is
a general locally compact, separable metric space, we discuss
properties of the model
\begin{eqnarray*}
\left\{
\begin{array}{l}
{\cal E}(F,G)=\frac{1}{2}\int_{{\cal P}_1(S)}\langle \nabla
F(\mu),\nabla G(\mu)\rangle_{\mu}
\Pi_{\alpha,\theta,\nu_0}(d\mu),\ \ F,G\in
{\cal F},\\
{\cal F}=\{F(\mu)=f(\langle \phi_1,\mu\rangle,\dots,\langle
\phi_d,\mu\rangle): \phi_i\in B_b(S),1\le i\le d,f\in
C^{\infty}(\mathbb{R}^d),d\ge 1\}.
\end{array}
\right.
\end{eqnarray*}
In particular, we prove the Mosco convergence of its projection
forms.

\smallskip


\section{Introduction}

For any $0\leq \alpha <1$ and $\theta>-\alpha$, let $U_k$,
      $k=1,2,\dots$, be a sequence of independent
      random variables such that $U_k$ has $Beta(1-\alpha,\theta+ k\alpha)$ distribution.
      Set
      $$
      V^{\alpha,\theta}_1 = U_1,\  V^{\alpha,\theta}_n = (1-U_1)\cdots (1-U_{n-1})U_n,\  n \geq 2,
      $$
and let ${\bf P}(\alpha,\theta)=(P_1(\alpha,\theta),
P_2(\alpha,\theta),\dots)$ denote $(V^{\alpha,\theta}_1,
V^{\alpha,\theta}_2,\dots)$
      in descending order. The distribution of $(V_1^{\alpha,\theta}, V_2^{\alpha,\theta}, \dots)$ is called
      the two-parameter GEM distribution, denoted by $GEM(\alpha,\theta)$. The law of ${\bf P}(\alpha,\theta)$ is called the two-parameter Poisson-Dirichlet
      distribution, denoted by $PD(\alpha,\theta)$ (\cite{PY}). For a locally compact, separable metric space $S$, and a sequence of i.i.d. $S$-valued random variables $\xi_k$,
      $k=1,2,\dots$, with common
       distribution $\nu_0$ on $S$, let
      $$
      \Theta_{\alpha,\theta, \nu_0}=\sum_{k=1}^{\infty}P_k(\alpha,\theta)\delta_{\xi_k}.
      $$
      Hereafter, we denote by $\delta_x$ the Dirac delta measure at $x$ for $x\in S$. The distribution of $\Theta_{\alpha,\theta, \nu_0}$, denoted by $Dirichlet(\alpha,\theta,\nu_0)$ or
      $\Pi_{\alpha,\theta,\nu_0}$, is called the two-parameter Dirichlet process.  Both $GEM(\alpha,\theta)$ and  $PD(\alpha,\theta)$  carry the information on proportions only while $\Pi_{\alpha,\theta,\nu_0}$ contains information on both proportions and types or labels.

      The two-parameter models are natural generalizations to the case $\alpha=0$.
      Specifically $PD(0,\theta)$, $GEM(0,\theta)$ and $\Pi_{0, \theta,\nu_0}$ correspond to the well known Poisson-Dirichlet distribution,
      the GEM distribution and the Dirichlet process,
      respectively.       The Poisson-Dirichlet distribution $PD(0,\theta)$ was introduced by Kingman in \cite{Kingman75}
 to describe the distribution of gene frequencies in a large neutral
      population at a particular locus. The component $P_k(\theta)$ represents the
      proportion of the $k$-th most frequent allele. The age-ordered proportions follow the GEM distribution.
      The Dirichlet process  $\Pi_{0,\theta,\nu_0}$ first appeared in \cite{Fer73} in the context of Bayesian statistics.
      It is a pure atomic random measure with masses distributed according to
      $PD(0,\theta)$. In the context of population genetics, both the Poisson-Dirichlet distribution and the Dirichlet process appear
   as approximations to the equilibrium behavior of certain large populations evolving under the influence of mutation
   and random genetic drift.

Let
\[
\nabla_{\infty}:=\left\{(x_1,x_2,\ldots): x_1\geq x_2\geq \cdots
\geq 0, \sum_{i=1}^{\infty}x_i=1 \right\}
\] denote the infinite
dimensional ordered simplex and
$${\overline{\nabla}}_{\infty}:=\left\{(x_1,x_2,\dots):x_1\ge
x_2\ge\cdots\ge 0,\ \sum_{i=1}^{\infty}x_i\le 1\right\}$$ be the
closure of ${\nabla}_{\infty}$ in the product space
$[0,1]^{\infty}$. In \cite{EK81} an infinite dimensional diffusion
process, the unlabeled infinitely-many-neutral-alleles model, is
constructed on $\overline{\nabla}_{\infty}$ with generator
$$
{\cal A}_{\theta}
=\frac{1}{2}\left\{\sum_{i,j=1}^{\infty}x_i(\delta_{ij}-x_j)\frac{\partial^2}{\partial
x_i\partial x_j} -\sum_{i=1}^{\infty}\theta x_i \frac{\partial
}{\partial x_i}\right\},
$$
defined on an appropriate domain. The reversible measure of this
process is shown to be $PD(0,\theta)$.

Let $d\ge 1$, $\phi_i\in B_b(S)$, $1\le i\le d$, $f\in
C^{\infty}(\mathbb{R}^d)$ and $F(\mu)=f(\langle
\phi_1,\mu\rangle,\dots,\langle \phi_d,\mu\rangle)$ for $\mu\in
{\cal P}_1(S)$. Hereafter, we denote by $B_b(S)$ the set of all
bounded Borel measurable functions on  $S$,
$C^{\infty}(\mathbb{R}^d)$ the set of all infinitely
differentiable functions on $\mathbb{R}^d$, and ${\cal P}_1(S)$
the space of all probability measures on the Borel
$\sigma$-algebra ${\cal B}(S)$ in $S$. For $x\in S$ and $\mu\in
{\cal P}_1(S)$, we define
\begin{eqnarray*}
\nabla_xF(\mu)&:=&\left.\frac{dF}{ds}(\mu+s\delta_x)\right|_{s=0}\\
&=&\sum_{i=1}^d\partial_if(\langle \phi_1,\mu\rangle,\dots,\langle
\phi_d,\mu\rangle)\phi_i(x).
\end{eqnarray*}
We write $\nabla F(\mu)$ for the function $x\rightarrow
\nabla_xF(\mu)$. For $\phi,\psi\in B_b(S)$, define
$$
\langle \phi,\psi\rangle_{\mu}:=\int_S \phi\psi d\mu-\left(\int_S
\phi d\mu\right)\left(\int_S \psi d\mu\right).
$$

 Given $\nu_0\in{\cal
P}_1(S)$ we consider the operator $A$ of the form
$$Ag(x) = \frac{\theta}{2}\int(g(y)-g(x))\nu_0(dy),\ \ g \in B_b(S).$$
Then, the Fleming-Viot process (cf. \cite{FleVio79} and
\cite{EtKu93}) with neutral parent independent mutation or the
labeled infinitely-many-neutral-alleles model is a pure atomic
measure-valued Markov process with generator
 \[
 L_{\theta} F(\mu)  =\frac{1}{2}\sum_{i,j=1}^d\partial_i\partial_j f(\langle
\phi_1,\mu\rangle,\dots,\langle
\phi_d,\mu\rangle)\langle\phi_i,\phi_j\rangle_{\mu}+ \langle
A\nabla F(\mu)(\cdot),\mu \rangle.\] For compact space $S$  and
diffuse probability $\nu_0$, i.e., $\nu_0(x)=0$ for every $x$ in
$S$, it is known (\cite{Et90}) that the labeled
infinitely-many-neutral-alleles model is time-reversible with
reversible measure $\Pi_{0,\theta,\nu_0}$.

 It is natural to ask whether these diffusion processes have two-parameter analogues when $\alpha$ is positive.
 Many progresses have been made in this direction over the last decade.  In \cite{fenwang07}, a class of infinite dimensional reversible diffusions is constructed and the reversible measure is $GEM(\alpha,\theta)$.
 The unlabeled infinitely-many-neutral-alleles model in \cite{EK81} is generalized to the two-parameter setting in \cite{P} where the generator of the process on appropriate domain has the form
$$
{\cal A}_{\alpha,\theta}
=\frac{1}{2}\left\{\sum_{i,j=1}^{\infty}x_i(\delta_{ij}-x_j)\frac{\partial^2}{\partial
x_i\partial x_j} -\sum_{i=1}^{\infty}(\alpha+\theta x_i) \frac{\partial
}{\partial x_i}\right\},
$$
and the reversible measure turns out to be $PD(\alpha, \theta)$.
The process, called Petrov diffusion, is derived as the continuum
limit of a family of up-down Markov  chains involving the Chinese
restaurant process. Connections to Bayesian statistics and ecology
are explored in \cite{Rugg14} and \cite{R-W09}. Going back to the
context of population genetics, the Petrov diffusion is
constructed recently in \cite{CBERS17} from a family of the
Wright-Fisher diffusions with special selection scheme.  In
\cite{FPRW17}, two interval partition-valued diffusions are
constructed and the corresponding stationary distributions are
$PD(1/2,0)$ and $PD(1/2,1/2)$, the two cases that are connected to the excursion intervals of Brownian motion and Brownian bridge (\cite{PPY92}, \cite{PY92}).

The situation is more complex in the construction of the labelled diffusion processes in the two-parameter setting. The only model we know of is the one in \cite{feng-sun10} where the type space consists of two types. In the case $\alpha=0$, the Dirichlet process $\Pi_{0,\theta,\nu_0}$ has the partition property, i.e., projection of $\Pi_{0,\theta,\nu_0}$ on any finite partition of the type space $S$ is a Dirichlet distribution. Exploring the connection between the Wright-Fisher diffusion and the Dirichlet distribution one can naturally construct the Fleming-Viot process from the finite-dimensional Wright-Fisher diffusions. When $\alpha$ is positive, the projection $\Pi_{\alpha,\theta,\nu_0}$ on any finite partition of $S$ has a complicated distribution in general, and finite dimensional diffusion models are no longer available.

  The main objective of this paper is to find a labelled reversible diffusion process with $\Pi_{\alpha,\theta,\nu_0}$ as
  the reversible measure for certain positive $\alpha$. This can be viewed as a two-parameter generalization of the Fleming-Viot process with parent independent mutation.
  The range of parameters we consider throughout the paper is $\alpha=1/2$ and $\theta>-1/2$.

In Section 2, we construct the process when the base measure
$\nu_0$ has countable support. Since the partition property does
not hold,  we will explore the partition structure through
Dirichlet forms. This allows us to avoid certain exceptional sets
that cause problems in the representation of generators. In
Section 3, we consider the general type space with diffuse base
measure. We first show that cylindrical functions do not belong to
the domain of the pre-generator of the classical bilinear form.
To establish the closability, we consider the relaxation of the
bilinear form. The process is then constructed, when  $S$ is a
compact Polish space, by taking the Mosco limit.

\section{Dynamic model with atomic base distribution}\setcounter{equation}{0}
Throughout this section, let $S=\mathbb{N}$, the set of all natural numbers. We consider the bilinear form
\begin{eqnarray*}
\left\{
\begin{array}{l}
{\cal E}(F,G)=\frac{1}{2}\int_{{\cal P}_1(\mathbb{N})}\langle
\nabla F(\mu),\nabla G(\mu)\rangle_{\mu}
\Pi_{\alpha,\theta,\nu_0}(d\mu),\ \ F,G\in
{\cal F},\\
{\cal F}=\{F(\mu)=f(\mu(1),\dots,\mu(d)):f\in
C^{\infty}(\mathbb{R}^d), d\ge 1\}.
\end{array}
\right.
\end{eqnarray*}
\begin{thm}\label{thm90}
The bilinear form $({\cal E},{\cal F})$ is closable on $L^2({\cal
P}_1(\mathbb{N});\Pi_{\alpha,\theta,\nu_0})$ and its closure
$({\cal E}, D({\cal E}))$ is a quasi-regular Dirichlet form. The
 diffusion process associated with $({\cal E},
D({\cal E}))$ is time-reversible with the stationary distribution
$\Pi_{\alpha,\theta,\nu_0}$.
\end{thm}

Before proving Theorem \ref{thm90}, we make some preparation.

For $d\ge 1$, we define
$$
\Delta_d:=\left\{(x_1,\dots, x_d)\in \mathbb{R}^{d}:x_i\ge 0\ {\rm
and}\ \sum_{i=1}^dx_i\le 1\right\}.
$$
Denote $p_i=\nu_0(i)$, $1\le i\le d$, and
$p_{d+1}=1-p_1-\cdots-p_d$. For $(x_1,\dots, x_d)\in \Delta_d$, we
define
$$ \rho_d(x_1,\dots,x_d):=\frac{p_1\cdots
p_{d+1}\Gamma(\theta+\frac{d+1}{2})}{\pi^{d/2}\Gamma(\theta+\frac{1}{2})}\frac{x_1^{-3/2}\cdots
x_{d+1}^{-3/2}}{\left(\frac{p_1^2}{x_1}+\cdots+
\frac{p_{d+1}^2}{x_{d+1}}\right)^{\theta+\frac{d+1}{2}}}1_{\{x_1+\cdots
+x_{d+1}=1\}}.
$$

 Denote $S_{d+1}=\{(x_1,\dots,
x_{d+1})\in \mathbb{R}^{d+1}:x_i\ge 0\ {\rm and}\
\sum_{i=1}^{d+1}x_i= 1\}$. Following the argument of \cite[proof
of Lemma 3.1]{C}, we get
\begin{eqnarray*}
&&\int_{\Delta_d}\frac{\frac{p_1^4}{x_1^2}}{\left(\frac{p_1^2}{x_1}+\cdots+
\frac{p_d^2}{x_{d}}+\frac{p^2_{d+1}}{1-x_1-\cdots-x_{d}}\right)^2}\rho_d(x_1,\dots, x_d)dx_1\cdots dx_d\\
&=&\frac{p_1^4\Gamma(\theta+\frac{d+1}{2})p_1\cdots
p_{d+1}}{\pi^{d/2}\Gamma(\theta+\frac{1}{2})}\int_{S_{d+1}}\frac{x_1^{-7/2}x_2^{-3/2}\cdots
x_{d+1}^{-3/2}}{\left(\frac{p^2_1}{x_1}+\cdots+
\frac{p^2_{d+1}}{x_{d+1}}\right)^{2+\theta+\frac{d+1}{2}}}dx_1\cdots dx_{d+1}\\
&=&\frac{p_1^4\Gamma(\theta+\frac{d+1}{2})}{\pi^{d/2}\Gamma(\theta+\frac{1}{2})}\frac{p_1\cdots
p_{d+1}}{\Gamma(2+\theta+\frac{d+1}{2})2^{2+\theta+\frac{d+1}{2}}}\\
& &\cdot
\int_{0}^{\infty}\cdots \int_{0}^{\infty} (s_1+\cdots +s_{d+1})^{-\theta}s_1^{-7/2}\prod_{j>1}^{d+1}s_j^{-3/2}\prod_{j=1}^{d+1}e^{-p_j^2/2s_j}ds_1\cdots ds_{d+1}\\
&=&\frac{p_1^5\Gamma(\theta+\frac{d+1}{2})}{\pi^{d/2}\Gamma(\theta+\frac{1}{2})}\frac{(2\pi)^{d/2}}
{\Gamma(2+\theta+\frac{d+1}{2})2^{2+\theta+\frac{d+1}{2}}}\\
& &\cdot
\int_0^{\infty}\left[\int_0^{\infty}(s_1+u)^{-\theta}\frac{1-p_1}{\sqrt{2\pi}}e^{-\frac{(1-p_1)^2}{2u}}u^{-3/2}du\right]s_1^{-7/2}e^{-\frac{p_1^2}{2s_1}}ds_1\\
&\le&\frac{C(\theta,p_1)}{(d+1+2\theta)(d+3+2\theta)},
\end{eqnarray*}
where $C(\theta,p_1)$ is a positive constant depending only on
$\theta$ and $p_1$. Further, we obtain by symmetry that for each
$1\le i\le p$,
\begin{equation}\label{41}
\int_{\Delta_d}\frac{\frac{p_i^4}{x_i^2}}{\left(\frac{p_1^2}{x_1}+\cdots+
\frac{p_d^2}{x_{d}}+\frac{p^2_{d+1}}{1-x_1-\cdots-x_{d}}\right)^2}\rho_d(x_1,\dots,
x_d)dx_1\cdots
dx_d\le\frac{C(\theta,p_i)}{(d+1+2\theta)(d+3+2\theta)},
\end{equation}
where $C(\theta,p_i)$ is a positive constant depending only on
$\theta$ and $p_i$.

Denote $C^{\infty}({\Delta_d}):=\{f|_{\Delta_d}:f\in
C^{\infty}(\mathbb{R}^d)\}$ and  ${\Delta^\circ_d}:=$ the interior
of ${\Delta_d}$.  For $f\in C^{\infty}({\Delta_d})$, we define
\begin{eqnarray}\label{3}
L^{(d)}f(x)=\left\{
\begin{array}{ll}
\frac{1}{2}\sum\limits_{i=1}^dx_i\partial^2_if(x)-\frac{1}{2}\sum\limits_{i,j=1}^dx_ix_j\partial_i\partial_jf(x)\\
\ \ \ \ +\frac{1}{2}\sum\limits_{i=1}^d\left[-\frac{1}{2}-\theta
x_i+\frac{(\theta+\frac{d+1}{2})\frac{p_i^2}{x_i}}{\frac{p_1^2}{x_1}+\cdots+
\frac{p_{d}^2}{x_{d}}+\frac{p_{d+1}^2}{1-x_1-\cdots-x_{d}}}\right]\partial_if(x),&\ \ x\in {\Delta^\circ_d},\\
0,&\ \ x\in \Delta_d\backslash{\Delta^\circ_d}.
\end{array}
\right.
\end{eqnarray}

If $f\in C^{\infty}(\mathbb{R}^p)$ for some $p\le d$, we regard
$f$ as a function in $C^{\infty}(\mathbb{R}^d)$ by setting
$f(x)=f(x_1,\dots,x_p)$ for $x=(x_1,\dots,x_d)\in \mathbb{R}^d$. By
(\ref{3}), there exists a constant $C(p,f)>0$, which
 depends on $p$, $f$ and is independent of $d$, such that for
any $x\in {\Delta^\circ_d}$,
\begin{equation}\label{42}
|L^{(d)}f(x)|\le
C(p,f)\left[1+(d+1)\sum_{i=1}^p\frac{\frac{p_i^2}{x_i}}{\frac{p_1^2}{x_1}+\cdots+
\frac{p_d^2}{x_{d}}+\frac{p^2_{d+1}}{1-x_1-\cdots-x_{d}}}\right].
\end{equation}
 By (\ref{41}) and (\ref{42}), we get
\begin{equation}\label{43}\int_{\Delta_d}|L^{(d)}f(x)|^2\rho_d(x) dx\le
C^*(\theta,\nu_0(1),\dots,\nu_0(p),p,f),
\end{equation}
where $C^*(\theta,\nu_0(1),\dots,\nu_0(p),p,f)$ is a positive
constant depending only on $\theta,\nu_0(1),\dots,\nu_0(p),p$, and
$f$.
\vskip 0.5cm
\noindent {\bf Proof of Theorem \ref{thm90}.}\quad Let $d\ge 1$. We consider the map
\begin{eqnarray*}
\Upsilon_d:\hskip -0.4cm&&{\cal
P}_1(\mathbb{N})\longrightarrow\Delta_d,\nonumber\\
&&\mu\longrightarrow\Upsilon_d(\mu)=(\mu(1),\dots,\mu(d)).
\end{eqnarray*}
By \cite[Theorem 3.1]{C}, we have
$\Pi_{\alpha,\theta,\nu_0}\circ\Upsilon_d^{-1}=\rho_d(x_1,\dots,x_d)dx_1\cdots
dx_d$.
 The induced bilinear form of $({\cal E},{\cal F})$ by the map
$\Upsilon_d$ is given by
\begin{equation}\label{cale}
{\cal E}^{(d)}(f,g)=\frac{1}{2}\sum_{i,j=1}^d\int_{\Delta_d}
x_i(\delta_{ij}-x_j)\partial_if(x)\partial_jg(x)\rho_d(x) dx,\ \
f,g\in C^{\infty}({\Delta_d}).
\end{equation}

For $x=(x_1,\dots,
x_d)\in \Delta_d$ and $1\le j\le d$, we define
$$
V_j(x)=\left(\sum_{i=1}^dx_i(\delta_{ij}-x_j)\partial_if(x)\right)\rho_d(x)g(x),
$$
and $$V=(V_1,\dots,V_d).$$ Denote by $\partial\Delta_d$ the
boundary of $\Delta_d$, $\mathbf{n}$ the outward pointing unit
normal field of $\partial\Delta_d$, and $d{\cal S}_d$ the induced
volume form on the surface $\partial\Delta_d$. For the face
$\{x=(x_1,\dots, x_d)\in \Delta_d:x_j=0\}$, $1\le j\le d$, we have
\begin{eqnarray*}
V\cdot \mathbf{n}&=&V_j\\
&=&\left(x_j(1-x_j)\partial_jf(x)-\sum_{i\not=
j}x_ix_j\partial_if(x)\right)\rho_d(x)g(x)\\
&=&0,
\end{eqnarray*}
and for the face $\{x=(x_1,\dots, x_d)\in
\Delta_d:\sum_{j=1}^dx_j=1\}$, we have
\begin{eqnarray*}
V\cdot \mathbf{n}&=&\frac{1}{\sqrt{d}}\sum_{j=1}^dV_j\\
&=&\frac{1}{\sqrt{d}}\left(\sum_{i=1}^dx_i\partial_i f(x)\right)\left(1-\sum_{j=1}^dx_j\right)\rho_d(x)g(x)\\
&=&0.
\end{eqnarray*}
Hence $V\cdot \mathbf{n}=0$ on $\partial\Delta_d$. Then, we obtain
by (\ref{43}) and the divergence theorem that
\begin{eqnarray*}
&&{\cal E}^{(d)}(f,g)+\int_{\Delta_d}L^{(d)}f(x)g(x)\rho_d(x) dx\nonumber\\
&=&\int_{\Delta_d}\sum_{j=1}^d\partial_{j}\left[\left(\sum_{i=1}^dx_i(\delta_{ij}-x_j)\partial_if(x_1,\dots,
x_d)\right)\rho_d(x_1,\dots, x_d)g(x_1,\dots,
x_d)\right]dx_1\cdots dx_d\nonumber\\
&=&\int_{\Delta_d}{\rm div}Vdx_1\cdots dx_d\nonumber\\
&=&\int_{\partial\Delta_d}V\cdot \mathbf{n}d{\cal S}_d\nonumber\\
&=&0.
\end{eqnarray*}
Thus, we have
\begin{equation}\label{clo}
{\cal E}^{(d)}(f,g)=\int_{\Delta_d}-L^{(d)}f(x)g(x)\rho_d(x) dx.
\end{equation}

Now we use the estimate (\ref{43}) to show that
$({\cal E},{\cal F})$ is closable on $L^2({\cal
P}_1(\mathbb{N});\Pi_{\alpha,\theta,\nu_0})$. To this end, let
$\{F_n\in {\cal F}\}$ be a sequence satisfying
$$
\lim_{n\rightarrow\infty}\|F_n\|_{L^2({\cal
P}_1(\mathbb{N});\Pi_{\alpha,\theta,\nu_0})}=0\ {\rm and}\
\lim_{n,m\rightarrow\infty}{\cal E}(F_n-F_m,F_n-F_m)=0.
$$
Note that
\begin{eqnarray*}
{\cal E}(F_n,F_n)&=&{\cal E}(F_n-F_k,F_n)+{\cal E}(F_k,F_n)\\
&\le&{\cal E}^{1/2}(F_n-F_k,F_n-F_k){\cal E}^{1/2}(F_n,F_n)+{\cal
E}(F_k,F_n).
\end{eqnarray*}
To show $\lim_{n\rightarrow\infty}{\cal E}(F_n,F_n)=0$, we need
only show that for any fixed $k$,
$$
\lim_{n\rightarrow\infty}{\cal E}(F_k,F_n)=0.
$$

Suppose that $F_n(\mu)=f^{(n)}(\mu(1),\dots,\mu(p^{(n)}))$ with
$f^{(n)}\in C^{\infty}(\mathbb{R}^{p^{(n)}})$ and
$p^{(n)},n\in\mathbb{N}$. By (\ref{43}) and (\ref{clo}), we get
\begin{eqnarray*}
|{\cal E}(F_k,F_n)|&=&\left|{\cal E}^{(p^{(k)}\vee p^{(n)})}(f^{(k)},f^{(n)})\right|\\
&=&\left|(-L^{(p^{(k)}\vee p^{(n)})}f^{(k)},f^{(n)})_{L^2\left({\Delta_{(p^{(k)}\vee p^{(n)})}};\,\Pi_{\alpha,\theta,\nu_0}\circ\Upsilon^{-1}_{(p^{(k)}\vee p^{(n)})}\right)}\right|\\
&\le&
\sqrt{C^*(\theta,\nu_0(1),\dots,\nu_0(p^{(k)}),p^{(k)},f^{(k)})}\cdot\|F_n\|_{L^2({\cal
P}_1(\mathbb{N});\Pi_{\alpha,\theta,\nu_0})}\\
&\rightarrow&0\ \ {\rm as}\ n\rightarrow\infty.
\end{eqnarray*}
Thus, $({\cal E},{\cal F})$ is closable on $L^2({\cal
P}_1(\mathbb{N});\Pi_{\alpha,\theta,\nu_0})$.

Following the argument of \cite[Proposition 5.11 and Lemma
7.5]{S}, we can show that the closure $({\cal E}, D({\cal E}))$ of
$({\cal E},{\cal F})$ is a quasi-regular, symmetric, local
Dirichlet form on $L^2({\cal P}_1(\mathbb{N});
\Pi_{\alpha,\theta,\nu_0})$. Therefore, there exists an associated
diffusion process in ${\cal P}_1(\mathbb{N})$ which is
time-reversible with the stationary distribution
$\Pi_{\alpha,\theta,\nu_0}$.\hfill\fbox

Denote by $(L,D(L))$ the generator of $({\cal E}, D({\cal E}))$ on
$L^2({\cal P}_1(\mathbb{N}); \Pi_{\alpha,\theta,\nu_0})$. In the
following, we will give an explicit expression for $L$.
\begin{thm} {\rm (i)} ${\cal F}\subset D(L)$.

{\rm (ii)} For each $i\in \mathbb{N}$,
\begin{equation}\label{red}
\lim_{d\rightarrow\infty}\frac{(d+1)\frac{\nu_0(i)^2}{\mu(i)}}{\frac{\nu_0(1)^2}{\mu(1)}+\cdots+
\frac{\nu_0(d)^2}{\mu(d)}+\frac{\nu_0(d+1)^2}{1-\mu(1)-\cdots-\mu(d)}}\
\ {\it exists\ in}\ L^2({\cal P}_1(\mathbb{N});
\Pi_{\alpha,\theta,\nu_0}).
\end{equation}

{\rm (iii)} For $i\in \mathbb{N}$, denote by $B_i(\mu)$ the
$L^2$-limit given in (\ref{red}). Let
$F(\mu)=f(\mu(1),\dots,\mu(d))$ with $f\in
C^{\infty}(\mathbb{R}^d)$ and $d\ge 1$. We have
\begin{eqnarray}\label{black}
LF(\mu)&=&\frac{1}{2}\sum_{i=1}^d\mu(i)\partial^2_if(\mu(1),\dots,\mu(d))-\frac{1}{2}\sum_{i,j=1}^d\mu(i)\mu(j)\partial_i\partial_jf(\mu(1),\dots,\mu(d))\nonumber\\
&&+\frac{1}{2}\sum_{i=1}^d\left[-\frac{1}{2}-\theta
\mu(i)+\frac{1}{2}B_i(\mu)\right]\partial_if(\mu(1),\dots,\mu(d)).
\end{eqnarray}
\end{thm}

\noindent {\bf Proof.}\quad
Let $F(\mu)=f(\mu(1),\dots,\mu(d))$ for some $f\in
C^{\infty}(\mathbb{R}^d)$ and $d\ge 1$. For
$G(\mu)=g(\mu(1),\dots,\mu(d'))$ with $g\in
C^{\infty}(\mathbb{R}^{d'})$ and $d'\ge 1$, we obtain by (\ref{43}) that
$$
|{\cal E}(F,G)|=|{\cal E}^{(d\vee d')}(f,g)|\le
(C^*(\theta,\nu_0(1),\dots,\nu_0(d),d,f))^{1/2}\|G\|_{L^2({\cal
P}_1(\mathbb{N});\Pi_{\alpha,\theta,\nu_0})}.
$$
Since $G$ is arbitrary, we conclude that $F\in D(L)$ by \cite[Chapter I, Proposition
2.16]{MR}.

For $n> d$, we regard $f$ as a function in
$C^{\infty}(\mathbb{R}^n)$ by setting $f(x)=f(x_1,\dots,x_d)$ for
$x=(x_1,\dots,x_n)\in \mathbb{R}^n$. We claim that
\begin{equation}\label{sdf}
LF=\lim_{n\rightarrow\infty}(L^{(n)}f)\circ\Upsilon_n\ \ {\rm in}\
L^2({\cal P}_1(\mathbb{N});\Pi_{\alpha,\theta,\nu_0}).
\end{equation}
In fact, it is easy to see that
$$
P_n(LF)=(L^{(n)}f)\circ\Upsilon_n\ \ {\rm for}\ n\ge d,
$$
where $P_n$ is the orthogonal projection of $L^2({\cal
P}_1(\mathbb{N});\Pi_{\alpha,\theta,\nu_0})$ onto the closure of
$\{G(\mu)=g(\mu(1),\dots,\mu(n)): g\in
C^{\infty}(\mathbb{R}^n)\}$. Since ${\cal F}$ is dense in
$L^2({\cal P}_1(\mathbb{N});\Pi_{\alpha,\theta,\nu_0})$, we obtain
(\ref{sdf}).

For $i\in \mathbb{N}$, let $F_i(\mu)=\mu(i)$ for $\mu\in {\cal
P}_1(\mathbb{N})$. By (\ref{sdf}), we get
$$
LF_i(\mu)=\frac{1}{2}\lim_{n\rightarrow\infty}
\frac{(\theta+\frac{n+1}{2})\frac{\nu_0(i)^2}{\mu(i)}}{\frac{\nu_0(1)^2}{\mu(1)}+\cdots+
\frac{\nu_0(n)^2}{\mu(n)}+\frac{\nu_0(n+1)^2}{1-\mu(1)-\cdots-\mu(n)}}.
$$
Hence, (\ref{red}) holds and $B_i=4LF_i$. Therefore, we obtain
(\ref{black}) by (\ref{3}) and (\ref{sdf}).\hfill\fbox

\begin{rem} For $f\in C^{\infty}(\mathbb{R})$ and $x\in (0,1)$, we have
$$
L^{(1)}f(x)=\frac{1}{2}x(1-x)f''(x)+\frac{1}{2}\left\{\theta+\frac{1}{2}-(2\theta+1)x\right\}f'(x).
$$
The eigenvalues of $L^{(1)}$ are $-i(i+2\theta)/2$ with
multiplicity 1, $i\in \mathbb{N}$. It deserves further
investigation to   characterize  the eigenvalues  of $L^{(d)}$ for
$d>1$.
\end{rem}

\section{Dynamic model with diffuse base distribution}\setcounter{equation}{0}

In this section, let $S$ be a general locally compact, separable
metric space and $\nu_0$ a diffuse probability measure on $S$. We
consider the classical bilinear form
\begin{eqnarray}\label{discussion}
\left\{
\begin{array}{l}
{\cal E}(F,G)=\frac{1}{2}\int_{{\cal P}_1(S)}\langle \nabla
F(\mu),\nabla G(\mu)\rangle_{\mu}
\Pi_{\alpha,\theta,\nu_0}(d\mu),\ \ F,G\in
{\cal F},\\
{\cal F}=\{F(\mu)=f(\langle \phi_1,\mu\rangle,\dots,\langle
\phi_d,\mu\rangle): \phi_i\in B_b(S),1\le i\le d,f\in
C^{\infty}(\mathbb{R}^d),d\ge 1\}.
\end{array}
\right.
\end{eqnarray}
If $({\cal E},{\cal F})$ is closable on $L^2({\cal
P}_1(S);\Pi_{\alpha,\theta,\nu_0})$, then following the argument
of \cite[Proposition 5.11 and Lemma 7.5]{S}, we can show that the
closure of $({\cal E},{\cal F})$ is a quasi-regular, symmetric,
local Dirichlet form  on $L^2({\cal P}_1(S);
\Pi_{\alpha,\theta,\nu_0})$. Therefore, there exists an associated
diffusion process in ${\cal P}_1(S)$ which is time-reversible with
the stationary distribution $\Pi_{\alpha,\theta,\nu_0}$.

Up to now we still cannot prove that $({\cal E},{\cal F})$ is
closable on $L^2({\cal P}_1(S);\Pi_{\alpha,\theta,\nu_0})$. In the
following, we will discuss properties of the model
(\ref{discussion}). We fix a sequence
$\{(B^k_1,\dots,B^k_{2^k})\}_{k=1}^{\infty}$ of partitions of $S$
satisfying the following conditions:

(1) $\nu_0(B^k_j)=1/2^k$, $1\le j\le 2^k$, $k\in\mathbb{N}$.

(2) $B^k_j=B^{k+1}_{2j-1}\cup B^{k+1}_{2j}$, $1\le j\le 2^k$,
$k\in\mathbb{N}$.

\subsection{${\cal F}\not\subset  D(L)$}

Denote by $(L,D(L))$ the pre-generator of $({\cal E},{\cal F})$ on
$L^2({\cal P}_1(S);\Pi_{\alpha,\theta,\nu_0})$. A special feature
of the model (\ref{discussion}) is that ${\cal F}\not\subset
D(L)$. More precisely, we have the following result.
\begin{pro}
Suppose that $\theta=0$. Let
$H(\mu)=\langle1_{B^1_1},\mu\rangle$. There does not exist $LH\in
L^2(M_1(S); \Pi_{\alpha,\theta,\nu_0})$ such that
\begin{equation}\label{mai}
{\cal E}(H,G)=\int_{{\cal P}_1(S)}-LH(\mu)
G(\mu){\Pi_{\alpha,\theta,\nu_0}}(d\mu),\ \ \forall G\in {\cal F}.
\end{equation}
\end{pro}
{\bf Proof.}\quad Let $d\ge 1$ and $p\le d$. For
$x=(x_1,\dots,x_d)\in\mathbb{R}^d$, we define
$$f(x)=x_1+\cdots +x_p.$$ Set $p_i=1/(d+1)$, $1\le i\le d+1$. We define $L^{(d)}$ as in
(\ref{3}). Then, we have
\begin{equation}\label{sdf2}
L^{(d)}f(x)=\frac{1}{4}\left(-p+(d+1)\sum_{i=1}^p\frac{\frac{1}{x_i}}{\frac{1}{x_1}+\cdots+
\frac{1}{x_{d}}+\frac{1}{1-x_1-\cdots-x_{d}}}\right).
\end{equation}

Following the argument of \cite[proof of Lemma 3.1]{C}, we get
\begin{eqnarray*}
&&\int_{\Delta_d}\frac{\frac{1}{x_1^2}}{\left(\frac{1}{x_1}+\cdots+
\frac{1}{x_{d}}+\frac{1}{1-x_1-\cdots-x_{d}}\right)^2}\rho_d(x_1\dots x_d)dx_1\cdots dx_d\nonumber\\
&=&\frac{\Gamma(\frac{d+1}{2})}{\pi^{d/2}\Gamma(\frac{1}{2})}\int_{S_{d+1}}\frac{x_1^{-7/2}x_2^{-3/2}\cdots
x_{d+1}^{-3/2}}{\left(\frac{1}{x_1}+\cdots+
\frac{1}{x_{d+1}}\right)^{2+\frac{d+1}{2}}}dx_1\cdots dx_{d+1}\nonumber\\
&=&\frac{\Gamma(\frac{d+1}{2})}{(d+1)^4\pi^{d/2}\Gamma(\frac{1}{2})}\frac{p_1\cdots
p_{d+1}}{\Gamma(2+\frac{d+1}{2})2^{2+\frac{d+1}{2}}}
\int_{0}^{\infty}\cdots \int_{0}^{\infty} s_1^{-7/2}\prod_{j>1}^{d+1}s_j^{-3/2}\prod_{j=1}^{d+1}e^{-p_j^2/2s_j}ds_1\cdots ds_{d+1}\nonumber\\
&=&\frac{\Gamma(\frac{d+1}{2})}{(d+1)^4\pi^{d/2}\Gamma(\frac{1}{2})}\frac{(2\pi)^{d/2}}
{(d+1)\Gamma(2+\frac{d+1}{2})2^{2+\frac{d+1}{2}}}
\int_0^{\infty}s_1^{-7/2}e^{-\frac{1}{2(d+1)^2s_1}}ds_1\nonumber\\
&=&\frac{3}{(d+1)(d+3)},
\end{eqnarray*}
and
\begin{eqnarray*}
&&\int_{\Delta_d}\frac{\frac{1}{x_1x_2}}{\left(\frac{1}{x_1}+\cdots+
\frac{1}{x_{d}}+\frac{1}{1-x_1-\cdots-x_{d}}\right)^2}\rho_d(x_1\dots x_d)dx_1\cdots dx_d\nonumber\\
&=&\frac{\Gamma(\frac{d+1}{2})}{\pi^{d/2}\Gamma(\frac{1}{2})}\int_{S_{d+1}}\frac{x_1^{-5/2}x_2^{-5/2}x_3^{-3/2}\cdots
x_{d+1}^{-3/2}}{\left(\frac{1}{x_1}+\cdots+
\frac{1}{x_{d+1}}\right)^{2+\frac{d+1}{2}}}dx_1\cdots dx_{d+1}\nonumber\\
&=&\frac{1}{(d+1)(d+3)}.
\end{eqnarray*}
Further, we obtain by symmetry that
\begin{equation}\label{sdf3}
\int_{\Delta_d}\frac{\frac{1}{x_i^2}}{\left(\frac{1}{x_1}+\cdots+
\frac{1}{x_{d}}+\frac{1}{1-x_1-\cdots-x_{d}}\right)^2}\rho_d(x_1\dots x_d)dx_1\cdots dx_d=\frac{3}{(d+1)(d+3)},\ \ 1\le i\le d,
\end{equation}
and
\begin{equation}\label{sdf4}
\int_{\Delta_d}\frac{\frac{1}{x_ix_j}}{\left(\frac{1}{x_1}+\cdots+
\frac{1}{x_{d}}+\frac{1}{1-x_1-\cdots-x_{d}}\right)^2}\rho_d(x_1\dots x_d)dx_1\cdots dx_d=\frac{1}{(d+1)(d+3)},\ \ 1\le i<j\le d.
\end{equation}
By (\ref{sdf2})--(\ref{sdf4}), we get
\begin{eqnarray}\label{4}
&&16\int_{\Delta_d}|L^{(d)}f(x)|^2\rho_d(x) dx\nonumber\\
&=&\int_{\Delta_d}\left(-p+(d+1)\sum_{i=1}^p\frac{\frac{1}{x_i}}{\frac{1}{x_1}+\cdots+
\frac{1}{x_{d}}+\frac{1}{1-x_1-\cdots-x_{d}}}\right)^2\rho_d(x_1\dots x_d)dx_1\cdots dx_d\nonumber\\
&=&-p^2+(d+1)^2p\cdot\frac{3}{(d+1)(d+3)}+(d+1)^2(p^2-p)\cdot\frac{1}{(d+1)(d+3)}\nonumber\\
&=&\frac{2p(d+1-p)}{d+3}.
\end{eqnarray}

Suppose there exists $LH\in L^2(M_1(S);
\Pi_{\alpha,\theta,\nu_0})$ such that (\ref{mai}) holds. For $k\in
\mathbb{N}$, we define
\begin{equation}\label{dis2}
{\cal F}_k=\{F(\mu)=f(\langle
1_{B^k_{1}},\mu\rangle,\dots,\langle
1_{B^k_{2^k-1}},\mu\rangle):f\in
C^{\infty}(\mathbb{R}^{2^k-1})\},
\end{equation}
and
$$
f_k(x_1,\dots,x_{2^k-1})=x_1+\cdots +x_{2^{k-1}},\ \
x_i\in\mathbb{R},\ 1\le i\le 2^k-1.
$$
Let $G(\mu)=g(\langle 1_{B^k_{1}},\mu\rangle,\dots,\langle
1_{B^k_{2^k-1}},\mu\rangle)$ for some $g\in
C^{\infty}(\mathbb{R}^{2^k-1})$. By (\ref{mai}), we get
\begin{eqnarray}\label{halm}
\int_{{\cal P}_1(S)}-LH(\mu)
G(\mu){\Pi_{\alpha,\theta,\nu_0}}(d\mu)&=&{\cal
E}(H,G)\nonumber\\
&=&{\cal
E}^{(2^k-1)}(f_k,g)\nonumber\\
&=&\int_{\Delta^{2^k-1}}-L^{(2^k-1)}f_k(x)g(x)\rho_{2^k-1}(x)dx.
\end{eqnarray}
Since $G\in {\cal F}_k$ is arbitrary, we obtain by (\ref{4}) and
(\ref{halm}) that
\begin{eqnarray*}
\|LH\|^2_{L^2({\cal
P}_1(S);\Pi_{\alpha,\theta,\nu_0})}&\ge&\|L^{(2^k-1)}f_k\|^2_{L^2({\Delta_{2^k-1}};\,\Pi_{\alpha,\theta,\nu_0}\circ\Upsilon^{-1}_{2^k-1})}\\
&=&\frac{2^{k-1}(2^k-2^{k-1})}{8(2^k+2)}.
\end{eqnarray*}
Since $k\in \mathbb{N}$ is arbitrary, there is a contradiction.
Therefore, there does not exist $LH\in L^2(M_1(S);
\Pi_{\alpha,\theta,\nu_0})$ such that (\ref{mai})
holds.\hfill\fbox

\subsection{Mosco convergence of projection forms}

Since we do not know if $({\cal E},{\cal F})$ is closable on
$L^2({\cal P}_1(S);\Pi_{\alpha,\theta,\nu_0})$, we consider the
relaxation of $({\cal E},{\cal F})$. By \cite[page 373]{M}, there
exists a greatest lower semicontinuous bilinear form on $L^2({\cal
P}_1(S);\Pi_{\alpha,\theta,\nu_0})$ which is a minorant of $({\cal
E},{\cal F})$. This unique determined closed form is called the
relaxation of $({\cal E},{\cal F})$, denoted by $(\Xi, D(\Xi))$.
We have that ${\cal F}\subset D(\Xi)$ and $\Xi(F,F)\le {\cal
E}(F,F)$ for any $F\in {\cal F}$, and for every $F\in D(\Xi)$,
\begin{eqnarray}\label{xcv}
\Xi(F,F)&=&\min\left\{\liminf_{n\rightarrow\infty}{\cal
E}(F_n,F_n):
F_n\in {\cal F}\ {\rm for}\ n\in \mathbb{N}\right.\nonumber\\
& &\ \ \ \ \ \ \ \ \ \ \ \ \ \ \ \left.{\rm and}\
\lim_{n\rightarrow\infty}F_n=F\ {\rm in}\ L^2({\cal
P}_1(S);\Pi_{\alpha,\theta,\nu_0})\right\}.
\end{eqnarray}
Note that if $({\cal E},{\cal F})$ is closable, then $(\Xi,
D(\Xi))$ is just the closure of $({\cal E},{\cal F})$ on
$L^2({\cal P}_1(S);\Pi_{\alpha,\theta,\nu_0})$.

By \cite[Corollary 2.8.2]{M}, $(\Xi, D(\Xi))$ is a Dirichlet form
on $L^2({\cal P}_1(S);\Pi_{\alpha,\theta,\nu_0})$. Further, if $S$
is a compact Polish space, then $(\Xi, D(\Xi))$ is a regular
Dirichlet form. Hence $(\Xi, D(\Xi))$ is associated with a Markov
process in ${\cal P}_1(S)$ which is time-reversible with the
stationary distribution $\Pi_{\alpha,\theta,\nu_0}$. Let ${\cal
F}_k$ be defined as in (\ref{dis2}) for $k\in \mathbb{N}$. In this
subsection, we will show that the limit of the sequence $\{({\cal
E},{\cal F}_k)\}$ is given by $(\Xi, D(\Xi))$.

For $k\in \mathbb{N}$, we consider the map
\begin{eqnarray*}
\Gamma_k:\hskip -0.4cm&&{\cal
P}_1(S)\longrightarrow\Delta_{2^k-1},\nonumber\\
&&\mu\longrightarrow\Gamma_k(\mu)=(\mu(B^k_1),\dots,\mu(B^k_{2^k-1})).
\end{eqnarray*}
For $F(\mu)=f(\langle 1_{B^k_{1}},\mu\rangle,\dots,\langle
1_{B^k_{2^k-1}},\mu\rangle)$ with $f\in
C^{\infty}(\mathbb{R}^{2^k-1})$ and $G(\mu)=g(\langle
1_{B^k_{1}},\mu\rangle,\dots,\langle 1_{B^k_{2^k-1}},\mu\rangle)$
with $g\in C^{\infty}(\mathbb{R}^{2^k-1})$, we obtain by
(\ref{clo}) that
\begin{eqnarray*}
{\cal E}(F,G)&=&{\cal
E}^{(2^k-1)}(f,g)\\
&=&\int_{\Delta_{2^k-1}}-L^{(2^k-1)}f(x)g(x)\rho_{2^k-1}(x) dx\\
&=&((-L^{(2^k-1)}f)\circ\Gamma_k,G)_{L^2({\cal
P}_1(S);\Pi_{\alpha,\theta,\nu_0})}.
\end{eqnarray*}
Hence $({\cal E},{\cal F}_k)$ is closable on $L^2({\cal
P}_1(S);\Pi_{\alpha,\theta,\nu_0})$ by \cite[Chapter I,
Proposition 3.3]{MR}. Denote by $({\cal E},D({\cal E})_k)$ the
closure of $({\cal E},{\cal F}_k)$. We have $D({\cal E})_1\subset
D({\cal E})_2\subset\cdots\subset L^2({\cal
P}_1(S);\Pi_{\alpha,\theta,\nu_0})$ and $({\cal E},D({\cal
E})_{k+1})$ is an extension of $({\cal E},D({\cal E})_k)$ for each
$k\in \mathbb{N}$.

For $k\in \mathbb{N}$, we define the resolvent $(G^k_{\beta})
_{\beta>0}$ of $({\cal E},D({\cal E})_k)$ by
\begin{equation}\label{res}
{\cal E}_{\beta}(G^k_{\beta}F,G)=(F,G)_{L^2({\cal
P}_1(S);\Pi_{\alpha,\theta,\nu_0})},\ \ \forall G\in D({\cal
E})_k,
\end{equation}
where ${\cal E}_{\beta}(F,G):={\cal E}(F,G)+\beta(F,G)_{L^2({\cal
P}_1(S);\Pi_{\alpha,\theta,\nu_0})}$ for $F,G\in D({\cal E})_k$.
Given $F\in L^2({\cal P}_1(S);\Pi_{\alpha,\theta,\nu_0})$, the
existence and uniqueness of $G^k_{\beta}F\in D({\cal E})_k$
satisfying (\ref{res}) follows from the Riesz representation
theorem.

Denote by $(G_{\beta}) _{\beta>0}$ the strongly continuous
contraction resolvent associated with the Dirichlet form $(\Xi,
D(\Xi))$ on $L^2({\cal P}_1(S);\Pi_{\alpha,\theta,\nu_0})$. We
have the following characterization of
  $(G_{\beta}) _{\beta>0}$ by virtue of $(G^k_{\beta}) _{\beta>0}$.
\begin{thm}\label{pro00}
For every $\beta>0$, the sequence of resolvent operators
$\{G^k_{\beta}\}$ converges to $G_{\beta}$ in the strong operator
topology.
\end{thm}
{\bf Proof.}\quad  We first show that for any subsequence $\{k'\}$
of $\{k\}$, there exists a subsequence $\{k^{''}\}$ of $\{k'\}$
such that for every $\beta>0$ the sequence
$\{G^{k^{''}}_{\beta}\}$ converges to a resolvent operator.

For $k\in \mathbb{N}$, we define
\begin{eqnarray*}
{\cal F}^*_k&=&\{F(\mu)=f(\langle
1_{B^k_{1}},\mu\rangle,\dots,\langle
1_{B^k_{2^k-1}},\mu\rangle):\\
& &\ \ \ \ \ \ \ \ \ \ f\ {\rm is\ a\ polynonmial\ on}\
\mathbb{R}^{2^k-1}\ {\rm with\ rational\ coefficients}\}.
\end{eqnarray*}
Denote
\begin{equation}\label{bnm}
{\cal F}^*=\bigcup_{k=1}^{\infty}{\cal F}^*_k.
\end{equation} Let
$\mathbb{Q}_+$ be the set of all positive rational numbers. Note
that
\begin{equation}\label{conver}
\|G^{l}_{\beta}F\|_{L^2({\cal
P}_1(S);\Pi_{\alpha,\theta,\nu_0})}\le\frac{1}{\beta}\|F\|_{L^2({\cal
P}_1(S);\Pi_{\alpha,\theta,\nu_0})},\ \ \forall F\in{\cal
F}_k,l\ge k,\beta>0.
\end{equation}
By the diagonal argument, there exists a subsequence $\{k^{''}\}$
of $\{k'\}$ such that
$$
G^{k^{''}}_{\beta}F\ {\rm converges\ weakly\ in}\ L^2({\cal
P}_1(S);\Pi_{\alpha,\theta,\nu_0})\ {\rm as}\
k^{''}\rightarrow\infty,\ \ \forall F\in {\cal F}^*, \beta\in
\mathbb{Q}_+.
$$
We fix such a subsequence $\{k^{''}\}$ and define
$$G^*_{\beta}F:=w-\lim_{k^{''}\rightarrow\infty}G^{k^{''}}_{\beta}F\ {\rm in}\ L^2({\cal
P}_1(S);\Pi_{\alpha,\theta,\nu_0}),\ \ F\in {\cal F}^*, \beta\in
\mathbb{Q}_+.
$$

Let $k{''}\le l{''}$. For $F\in {\cal F}^*$ and $\beta\in
\mathbb{Q}_+$, we have
\begin{eqnarray}\label{cvb}
&&{\cal
E}_{\beta}(G^{k^{''}}_{\beta}F-G^{l^{''}}_{\beta}F,G^{k^{''}}_{\beta}F-G^{l^{''}}_{\beta}F)\nonumber\\
&=&{\cal
E}_{\beta}(G^{k^{''}}_{\beta}F-G^{l^{''}}_{\beta}F,G^{k^{''}}_{\beta}F)-{\cal
E}_{\beta}(G^{k^{''}}_{\beta}F-G^{l^{''}}_{\beta}F,G^{l^{''}}_{\beta}F)\nonumber\\
&=&\{(F,G^{k^{''}}_{\beta}F)-(F,G^{k^{''}}_{\beta}F)\}-(G^{k^{''}}_{\beta}F-G^{l^{''}}_{\beta}F,F)\nonumber\\
&=&(G^{l^{''}}_{\beta}F,F)-(G^{k^{''}}_{\beta}F,F).
\end{eqnarray}
Hence
$$
\lim_{k^{''},l^{''}\rightarrow\infty}{\cal
E}_{\beta}(G^{k^{''}}_{\beta}F-G^{l^{''}}_{\beta}F,G^{k^{''}}_{\beta}F-G^{l^{''}}_{\beta}F)=0.
$$
Thus,
\begin{equation}\label{conver2}
G^*_{\beta}F=\lim_{k^{''}\rightarrow\infty}G^{k^{''}}_{\beta}F\
{\rm in}\ L^2({\cal P}_1(S);\Pi_{\alpha,\theta,\nu_0}),\ \ F\in
{\cal F}^*, \beta\in \mathbb{Q}_+.
\end{equation}

By (\ref{conver}) and (\ref{conver2}), we get
\begin{equation}\label{conver3}
\|\beta G^*_{\beta}F\|_{L^2({\cal
P}_1(S);\Pi_{\alpha,\theta,\nu_0})}\le\|F\|_{L^2({\cal
P}_1(S);\Pi_{\alpha,\theta,\nu_0})},\ \ \forall F\in {\cal F}^*,
\beta\in \mathbb{Q}_+.
\end{equation}
 For every $\beta\in \mathbb{Q}_+$, by (\ref{conver3}), we can
extend $\beta G^*_{\beta}$ to a continuous contraction operator on
$L^2({\cal P}_1(S);\Pi_{\alpha,\theta,\nu_0})$. Further, by
(\ref{conver2}), the resolvent equations for
$\{G^{k^{''}}_{\beta}\}$, and the density of ${\cal F}^*$ in
$L^2({\cal P}_1(S);\Pi_{\alpha,\theta,\nu_0})$, we can obtain a
collection of continuous operators $\{G^*_{\beta}\}_{\beta>0}$ on
$L^2({\cal P}_1(S);\Pi_{\alpha,\theta,\nu_0})$ satisfying

(i) $\|\beta G^*_{\beta}\|_{L^2({\cal
P}_1(S);\Pi_{\alpha,\theta,\nu_0})}\le 1,\ \ \forall \beta>0$.

(ii)
$G^*_{\beta}F=\lim_{k^{''}\rightarrow\infty}G^{k^{''}}_{\beta}F,\
\ \forall F\in L^2({\cal P}_1(S);\Pi_{\alpha,\theta,\nu_0}),
\beta>0$.

(iii)
$G^*_{\beta}-G^*_{\gamma}=(\gamma-\beta)G^*_{\beta}G^*_{\gamma},\
\ \forall \beta,\gamma>0$.

Let $F\in {\cal F}^*$ and $\beta\in \mathbb{Q}_+$. By (\ref{cvb}),
we find that $\{(\beta G^{k^{''}}_{\beta}F,F)_{L^2({\cal
P}_1(S);\Pi_{\alpha,\theta,\nu_0})}\}_{k^{''}=1}^{\infty}$ is an increasing sequence.
Then, we obtain by (\ref{conver2}) that
\begin{eqnarray}\label{equal}
\liminf_{\beta\in \mathbb{Q}_+,\beta\rightarrow\infty}(\beta
G^*_{\beta}F,F)_{L^2({\cal
P}_1(S);\Pi_{\alpha,\theta,\nu_0})}&\ge&\limsup_{k{{''}}\rightarrow\infty}\lim_{\beta\in
\mathbb{Q}_+,\beta\rightarrow\infty}(\beta
G^{k^{''}}_{\beta}F,F)_{L^2({\cal
P}_1(S);\Pi_{\alpha,\theta,\nu_0})}\nonumber\\
&=&\|F\|_{L^2({\cal P}_1(S);\Pi_{\alpha,\theta,\nu_0})}.
\end{eqnarray}
By (i) and (\ref{equal}), we get
$$
\lim_{\beta\in \mathbb{Q}_+,\beta\rightarrow\infty}(\beta
G^*_{\beta}F,F)_{L^2({\cal
P}_1(S);\Pi_{\alpha,\theta,\nu_0})}=\|F\|_{L^2({\cal
P}_1(S);\Pi_{\alpha,\theta,\nu_0})}.
$$
Then, we have \begin{eqnarray*}&&\limsup_{\beta\in
\mathbb{Q}_+,\beta\rightarrow\infty}\|\beta
G^*_{\beta}F-F\|^2_{L^2({\cal P}_1(S);\Pi_{\alpha,\theta,\nu_0})}\\
&\le&\limsup_{\beta\in
\mathbb{Q}_+,\beta\rightarrow\infty}2[\|F\|_{L^2({\cal
P}_1(S);\Pi_{\alpha,\theta,\nu_0})}-(\beta
G^*_{\beta}F,F)_{L^2({\cal P}_1(S);\Pi_{\alpha,\theta,\nu_0})}]\\
&=&0.
\end{eqnarray*}
Thus, we obtain by (i), (iii), and the density of ${\cal F}^*$ in
$L^2({\cal P}_1(S);\Pi_{\alpha,\theta,\nu_0})$ that

(iv) $\lim_{\beta\rightarrow\infty}\|\beta
G^*_{\beta}F-F\|_{L^2({\cal
P}_1(S);\Pi_{\alpha,\theta,\nu_0})}=0,\ \ \forall F\in L^2({\cal
P}_1(S);\Pi_{\alpha,\theta,\nu_0})$.

By (i), (iii), and (iv), we know that $(G^*_{\beta})_{\beta>0}$ is
a strongly continuous contraction resolvent on $L^2({\cal
P}_1(S);\Pi_{\alpha,\theta,\nu_0})$ (cf. \cite[Chapter I,
Definition 1.4]{MR}). Then, there exists a unique symmetric
Dirichlet form $(\Lambda, D(\Lambda))$ on $L^2({\cal
P}_1(S);\Pi_{\alpha,\theta,\nu_0})$ such that its resolvent is
given by $(G^*_{\beta})_{\beta>0}$, i.e.,
$$
\Lambda(G^*_{\beta}F,G)+\beta(G^*_{\beta}F,G)=(F,G)_{L^2({\cal
P}_1(S);\Pi_{\alpha,\theta,\nu_0})},\ \ \forall F\in {L^2({\cal
P}_1(S);\Pi_{\alpha,\theta,\nu_0})}, G\in D(\Lambda).
$$
By (ii) and \cite[Theorem 2.4.1]{M}, we find that $({\cal
E},D({\cal E})_{k^{''}})$ converges to $(\Lambda, D(\Lambda))$ in
the sense of Mosco convergence as $k^{''}\rightarrow\infty$, i.e.,

(a) For every $F_{k^{''}}\in D({\cal E})_{k^{''}}$ converging
weakly to $F\in D(\Lambda)$ in $L^2({\cal
P}_1(S);\Pi_{\alpha,\theta,\nu_0})$,
$$\liminf_{k^{''}\rightarrow\infty}{\cal
E}(F_{k^{''}},F_{k^{''}})\ge \Lambda(F,F).$$

(b) For every $F\in D(\Lambda)$, there exists $F_{k^{''}}\in
D({\cal E})_{k^{''}}$ converging strongly to $F\in D(\Lambda)$ in
$L^2({\cal P}_1(S);\Pi_{\alpha,\theta,\nu_0})$, such that
$$
\limsup_{k^{''}\rightarrow\infty}{\cal
E}(F_{k^{''}},F_{k^{''}})\le \Lambda(F,F).
$$

By (a), we know that $(\Lambda, D(\Lambda))$ is a minorant of
$({\cal E},{\cal F})$. By (b) and (\ref{xcv}), we obtain that $
D(\Lambda)\subset D(\Xi)$ and $\Xi(F,F)\le \Lambda(F,F)$ for $F\in
D(\Lambda)$. Since $(\Xi, D(\Xi))$ is the greatest closed form on
$L^2({\cal P}_1(S);\Pi_{\alpha,\theta,\nu_0})$ which is a minorant
of $({\cal E},{\cal F})$, we get $(\Lambda, D(\Lambda))=(\Xi,
D(\Xi))$. Then, we obtain by (ii) that
$$
G_{\beta}F=G^*_{\beta}F=\lim_{k^{''}\rightarrow\infty}G^{k^{''}}_{\beta}F,\
\ \forall F\in L^2({\cal P}_1(S);\Pi_{\alpha,\theta,\nu_0}),
\beta>0.
$$
Since the subsequence $\{k'\}$ of $\{k\}$ is arbitrary, we get
\begin{equation}\label{resolv}
G_{\beta}F=\lim_{k\rightarrow\infty}G^{k}_{\beta}F,\ \ \forall
F\in L^2({\cal P}_1(S);\Pi_{\alpha,\theta,\nu_0}), \beta>0.
\end{equation}\hfill\fbox

As a direct consequence of Theorem \ref{pro00} and \cite[Theorem
2.4.1]{M}, we obtain the Mosco convergence of projection forms of
the model (\ref{discussion}).
\begin{cor}
The sequence of bilinear forms $({\cal E},{\cal F}_k)$ converges
to $(\Xi, D(\Xi))$ in the sense of Mosco convergence.
\end{cor}

For $k\in\mathbb{N}$, we define  the bilinear form ${\cal
E}^{(2^k-1)}$ as in (\ref{cale}). By (\ref{clo}), we know that
$({\cal E}^{(2^k-1)},C^{\infty}({\Delta_{2^k-1}}))$ is closable on
$L^2({\Delta_{2^k-1}},\Pi_{\alpha,\theta,\nu_0}\circ\Gamma_{k}^{-1})$.
Denote by $({\cal E}^{(2^k-1)},D({\cal E}^{(2^k-1)}))$ the closure
of $({\cal E}^{(2^k-1)},C^{\infty}({\Delta_{2^k-1}}))$,
$(T^{(2^k-1)}_t)_{t\ge 0}$ the semigroup associated with $({\cal
E}^{(2^k-1)},D({\cal E}^{(2^k-1)}))$  on
$L^2({\Delta_{2^k-1}},\Pi_{\alpha,\theta,\nu_0}\circ\Gamma_{k}^{-1})$,
and  ${Q}_k$ the orthogonal projection of $L^2({\cal
P}_1(S);\Pi_{\alpha,\theta,\nu_0})$ onto the closure of ${\cal
F}_k$. For $F\in L^2({\cal P}_1(S);\Pi_{\alpha,\theta,\nu_0})$ and
$t\ge 0$, we define
$$
T^k_tF=(T^{(2^k-1)}_t(({Q}_kF)\circ\Gamma_{k}^{-1}))\circ\Gamma_{k}.
$$
Then, $(T^k_t)_{t\ge 0}$ is the semigroup associated with the
bilinear form $({\cal E},D({\cal E})_k)$ on $L^2({\cal
P}_1(S);\Pi_{\alpha,\theta,\nu_0})$.

Denote by $(T_t)_{t\ge 0}$ the strongly continuous contraction
semigroup associated with the Dirichlet form $(\Xi, D(\Xi))$ on
$L^2({\cal P}_1(S);\Pi_{\alpha,\theta,\nu_0})$. We have the
following characterization of
  $(T_t)_{t\ge 0}$ by virtue of $(T^k_t)_{t\ge 0}$.

\begin{thm}\label{corsemi}
For every $t\ge 0$, the sequence of semigroup operators
$\{T^k_t\}$ converges to $T_t$ in the strong operator
topology.
\end{thm}

\noindent {\bf Proof.}\quad   Let $\{k'\}$ be a subsequence
of $\{k\}$. By the diagonal argument, there exists a subsequence $\{k^{''}\}$
of $\{k'\}$ such that
$$w-\lim_{k^{''}\rightarrow\infty}T^{k^{''}}_tF\ {\rm exists}\ {\rm in}\ L^2({\cal
P}_1(S);\Pi_{\alpha,\theta,\nu_0}),\ \ \forall F\in {\cal F}^*, t\in
\mathbb{Q}_+,
$$
where ${\cal F}^*$ is defined as in (\ref{bnm}). We define
$$
T^{\Delta}_tF:=w-\lim_{k^{''}\rightarrow\infty}T^{k^{''}}_tF\ {\rm in}\ L^2({\cal
P}_1(S);\Pi_{\alpha,\theta,\nu_0}),\ \ F\in {\cal F}^*, t\in
\mathbb{Q}_+.
$$
By the density of ${\cal F}^*$ in $L^2({\cal
P}_1(S);\Pi_{\alpha,\theta,\nu_0})$ and the contraction of the
semigroup operators $\{T^k_t\}$, we can extend
$(T^{\Delta}_t)_{t\in \mathbb{Q}_+}$ to a collection of
contraction linear operators on $L^2({\cal
P}_1(S);\Pi_{\alpha,\theta,\nu_0})$ such that
\begin{equation}\label{addx}
T^{\Delta}_tF=w-\lim_{k^{''}\rightarrow\infty}T^{k^{''}}_tF\ {\rm in}\ L^2({\cal
P}_1(S);\Pi_{\alpha,\theta,\nu_0}),\ \ \forall F\in L^2({\cal
P}_1(S);\Pi_{\alpha,\theta,\nu_0}), {t\in
\mathbb{Q}_+}.
\end{equation}

By (\ref{addx}), we find that
\begin{equation}\label{addx111}
t\rightarrow(T^{\Delta}_tF,F)_{L^2({\cal P}_1(S);\Pi_{\alpha,\theta,\nu_0})}\ {\rm is\ decreasing\ on}\ \mathbb{Q}_+,\ \ \forall F\in L^2({\cal P}_1(S);\Pi_{\alpha,\theta,\nu_0}).
\end{equation} Hence, for any $t\ge 0$ and $F,G\in L^2({\cal P}_1(S);\Pi_{\alpha,\theta,\nu_0})$, we have
\begin{eqnarray}\label{addx112}
&&\lim_{s\in \mathbb{Q}_+,s\downarrow t}(T^{\Delta}_sF,G)_{L^2({\cal P}_1(S);\Pi_{\alpha,\theta,\nu_0})}\nonumber\\
&&\ \ \ \ \ \ =\frac{1}{4}
\lim_{s\in \mathbb{Q}_+,s\downarrow t}\{(T^{\Delta}_t(F+G),(F+G))_{L^2({\cal P}_1(S);\Pi_{\alpha,\theta,\nu_0})}\nonumber\\
&&\ \ \ \ \ \ \ \ \ \ \ \ \ \ \ \ \ \ -(T^{\Delta}_t(F-G),(F-G))_{L^2({\cal P}_1(S);\Pi_{\alpha,\theta,\nu_0})}\}\ {\rm exists}.\ \ \ \
\end{eqnarray}
By (\ref{addx111}) and (\ref{addx112}), we know that
\begin{equation}\label{addx11}
T^*_tF:=w-\lim_{s\in \mathbb{Q}_+,s\downarrow t}T^{\Delta}_sF,\ \ \forall t\ge 0,F\in L^2({\cal
P}_1(S);\Pi_{\alpha,\theta,\nu_0})
\end{equation}
is well-defined. Moreover, by (\ref{addx11}), we can show that $(T_t^*)_{t\ge 0}$  is
a collection of contraction linear operators on $L^2({\cal
P}_1(S);\Pi_{\alpha,\theta,\nu_0})$.

By (\ref{addx111}) and (\ref{addx11}), we find that
$$
t\rightarrow(T^*_tF,F)_{L^2({\cal P}_1(S);\Pi_{\alpha,\theta,\nu_0})}
\ {\rm is\ decreasing\ on}\ [0,\infty),\ \ \forall F\in L^2({\cal P}_1(S);\Pi_{\alpha,\theta,\nu_0}).
$$
Hence, there exists a collection of subsets $\{E_F\}_{F\in
L^2({\cal P}_1(S);\Pi_{\alpha,\theta,\nu_0})}$ of $[0,\infty)$
such that
$$
t\rightarrow(T^*_tF,F)_{L^2({\cal P}_1(S);\Pi_{\alpha,\theta,\nu_0})}
\ {\rm is\ continuous\ on}\ [0,\infty)\backslash E_F,\ \ \forall F\in L^2({\cal P}_1(S);\Pi_{\alpha,\theta,\nu_0}).
$$
For $t\ge 0$, we obtain by (\ref{addx}) and (\ref{addx11}) that
\begin{eqnarray}\label{semi1}
(T^*_tF,F)_{L^2({\cal
P}_1(S);\Pi_{\alpha,\theta,\nu_0})}&=&\lim_{s\in \mathbb{Q}_+,s\downarrow t}(T^{\Delta}_sF,F)_{L^2({\cal
P}_1(S);\Pi_{\alpha,\theta,\nu_0})}\nonumber\\
&=&\lim_{s\in \mathbb{Q}_+,s\downarrow t}\lim_{k^{''}\rightarrow\infty}(T^{k^{''}}_sF,F)_{L^2({\cal
P}_1(S);\Pi_{\alpha,\theta,\nu_0})}\nonumber\\
&\le&\lim_{k^{''}\rightarrow\infty}(T^{k^{''}}_tF,F)_{L^2({\cal
P}_1(S);\Pi_{\alpha,\theta,\nu_0})}.
\end{eqnarray}
For $t\in(0,\infty)\backslash E_F$ and $\varepsilon>0$, there
exists $\delta>0$ such that
\begin{equation}\label{delta1}(T^*_tF,F)_{L^2({\cal
P}_1(S);\Pi_{\alpha,\theta,\nu_0})}\ge (T^*_sF,F)_{L^2({\cal
P}_1(S);\Pi_{\alpha,\theta,\nu_0})}-\varepsilon,\ \ \forall s\in
((t-\delta)\vee0, t).\end{equation} By (\ref{addx}),
(\ref{addx11}), and (\ref{delta1}), we know that there exists
$t^*\in (0,t)\cap \mathbb{Q}_+$ such that
\begin{eqnarray*}
(T^*_tF,F)_{L^2({\cal
P}_1(S);\Pi_{\alpha,\theta,\nu_0})}&\ge&\lim_{k^{''}\rightarrow\infty}(T^{k^{''}}_{t^*}F,F)_{L^2({\cal
P}_1(S);\Pi_{\alpha,\theta,\nu_0})}-2\varepsilon\nonumber\\
&\ge&\lim_{k^{''}\rightarrow\infty}(T^{k^{''}}_tF,F)_{L^2({\cal
P}_1(S);\Pi_{\alpha,\theta,\nu_0})}-2\varepsilon.
\end{eqnarray*}
Since $\varepsilon$ is arbitrary, we get
\begin{eqnarray}\label{semi2}
(T^*_tF,F)_{L^2({\cal P}_1(S);\Pi_{\alpha,\theta,\nu_0})}\ge
\lim_{k^{''}\rightarrow\infty}(T^{k^{''}}_tF,F)_{L^2({\cal
P}_1(S);\Pi_{\alpha,\theta,\nu_0})},\ \ \forall
t\in(0,\infty)\backslash E_F.
\end{eqnarray}
By (\ref{semi1}) and (\ref{semi2}),  we get
\begin{eqnarray}\label{addxxc}
&&(T^*_tF,F)_{L^2({\cal
P}_1(S);\Pi_{\alpha,\theta,\nu_0})}=\lim_{k^{''}\rightarrow\infty}(T^{k^{''}}_tF,F)_{L^2({\cal
P}_1(S);\Pi_{\alpha,\theta,\nu_0})},\nonumber\\
& &\ \ \ \ \ \ \ \ \ \ \ \ \ \ \ \ \ \ \ \ \ \ \ \ \ \ \ \ \forall
t\in(0,\infty)\backslash E_F,F\in L^2({\cal
P}_1(S);\Pi_{\alpha,\theta,\nu_0}).
\end{eqnarray}

For $\beta>0$ and $F\in L^2({\cal
P}_1(S);\Pi_{\alpha,\theta,\nu_0})$, we obtain by (\ref{resolv}),
(\ref{addxxc}), and the dominated convergence theorem that
\begin{eqnarray}\label{right}
(G_{\beta}F,F)_{L^2({\cal P}_1(S);\Pi_{\alpha,\theta,\nu_0})}
&=&\lim_{k^{''}\rightarrow\infty}(G^{k^{''}}_{\beta}F,F)_{L^2({\cal P}_1(S);\Pi_{\alpha,\theta,\nu_0})}\nonumber\\
&=&\lim_{k^{''}\rightarrow\infty}
\int_0^{\infty}e^{-\beta t}(T_t^{k^{''}}F,F)_{L^2({\cal P}_1(S);\Pi_{\alpha,\theta,\nu_0})}dt\nonumber\\
&=&\int_0^{\infty}e^{-\beta t}(T_t^*F,F)_{L^2({\cal P}_1(S);\Pi_{\alpha,\theta,\nu_0})}dt.
\end{eqnarray}
By (\ref{right}), the right continuity of the function
$t\rightarrow(T^*_tF,F)_{L^2({\cal
P}_1(S);\Pi_{\alpha,\theta,\nu_0})}$ on $[0,\infty)$, and the
uniqueness of the Laplace transform, we find that
$$
(T_tF,F)_{L^2({\cal P}_1(S);\Pi_{\alpha,\theta,\nu_0})}
=(T_t^*F,F)_{L^2({\cal P}_1(S);\Pi_{\alpha,\theta,\nu_0})},\ \ \forall t\ge 0,F\in L^2({\cal
P}_1(S);\Pi_{\alpha,\theta,\nu_0}),
$$
which implies that
\begin{equation}\label{king}
T_tF=T^*_tF,\ \ \forall t\ge 0,F\in L^2({\cal
P}_1(S);\Pi_{\alpha,\theta,\nu_0}).
\end{equation}

By (\ref{addxxc}), (\ref{king}), the fact that the function
$t\rightarrow(T^{k^{''}}_tF,F)_{L^2({\cal P}_1(S);\Pi_{\alpha,\theta,\nu_0})}$ is decreasing on $[0,\infty)$,
and the continuity of the function
$t\rightarrow(T_tF,F)_{L^2({\cal P}_1(S);\Pi_{\alpha,\theta,\nu_0})}$ on $[0,\infty)$, we get
$$(T_tF,F)_{L^2({\cal
P}_1(S);\Pi_{\alpha,\theta,\nu_0})}=\lim_{k^{''}\rightarrow\infty}(T^{k^{''}}_tF,F)_{L^2({\cal
P}_1(S);\Pi_{\alpha,\theta,\nu_0})},\ \ \forall t\ge 0,F\in L^2({\cal
P}_1(S);\Pi_{\alpha,\theta,\nu_0}),
$$
which implies that
$$
T_tF=w-\lim_{k\rightarrow\infty}T^{k^{''}}_tF\ {\rm in}\ L^2({\cal
P}_1(S);\Pi_{\alpha,\theta,\nu_0}),\ \ \forall t\ge 0, F\in
L^2({\cal P}_1(S);\Pi_{\alpha,\theta,\nu_0}).
$$
Further, we obtain by the semigroup property that
$$
T_tF=\lim_{k\rightarrow\infty}T^{k^{''}}_tF\ {\rm in}\ L^2({\cal
P}_1(S);\Pi_{\alpha,\theta,\nu_0}),\ \ \forall t\ge 0, F\in
L^2({\cal P}_1(S);\Pi_{\alpha,\theta,\nu_0}).
$$
Since the subsequence $\{k'\}$ of $\{k\}$ is arbitrary, we get
$$
T_tF=\lim_{k\rightarrow\infty}T^{k}_tF\ {\rm in}\ L^2({\cal
P}_1(S);\Pi_{\alpha,\theta,\nu_0}),\ \ \forall t\ge 0, F\in L^2({\cal
P}_1(S);\Pi_{\alpha,\theta,\nu_0}).
$$\hfill\fbox

\bigskip


\end{document}